\documentclass[a4paper,12pt]{article}

\usepackage{latexsym}
\usepackage{amssymb,amsmath}
\usepackage{color}

\usepackage{scrtime}

\usepackage{hyperref}
\font\corsivo=rsfs10 at 12pt 
\font\doppio=msbm10 at 12pt 
\font\scdoppio=msbm8 

\newcommand{\R}{\hbox{\doppio R}}

\newcommand{\rn}{{\hbox{\scdoppio R}^N}}
\newcommand{\G}{\hbox{\doppio G}}     
\newcommand{\hei}{\hbox{\doppio H}} 
\newcommand{\RN}{{\R}^N}
\newcommand{\RL}{{\R}^l}

\newcommand{\A}{\mbox{\corsivo A}}      

\newcommand{\C}{\mbox{\corsivo C}}      

\newcommand{\grh}{\nabla_{\!\!H}}               
\newcommand{\lh}{\Delta_H}              

\newcommand{\Cuno}{\mbox{\corsivo C}^{\,1}}

\newcommand{\gr}{\nabla}
\newcommand{\grl}{\nabla_{\!\!L}}               
\newcommand{\diver}{\mathrm{div}}       
\newcommand{\diverl}{\diver_{\!L}}      
\newcommand{\decl}{:=}                  

\newcommand{\be}{\begin{equation}}      
\newcommand{\ee}{\end{equation}}        
\newcommand{\bern}{\begin{eqnarray*}}   
\newcommand{\eern}{\end{eqnarray*}}     

\newcommand{\norm}[1]{\abs{\abs{#1}}}

\newcommand{\SPC}{\mbox{{\bf S}-$p$-{\bf C}}}
\newcommand{\WPC}{\mbox{{\bf W}-$p$-{\bf C}}}

\newcommand{\MPC}{\mbox{{\bf M}-$p$-{\bf C}}}

\newtheorem{theorem}{Theorem}[section]
\newtheorem{lemma}[theorem]{Lemma}
\newtheorem{corollary}[theorem]{Corollary}
\newtheorem{definition}[theorem]{Definition}
\newtheorem{example}[theorem]{Example}
\newtheorem{remark}[theorem]{Remark}

\newcommand{\bp}{\noindent{\bf Proof. }}
\newcommand{\ep}{\logend\medskip}
\newcommand{\logend}{\hspace*{\fill}$\Box$}
\newcommand{\abs}[1]{\left |#1\right |}
\newcommand{\mint}{{\int\!\!\!\!\!\!-}}

\title{A new critical curve for a class of\\
 quasilinear elliptic systems}

\author{Lorenzo D'Ambrosio
 \\
{\small  Dipartimento di Matematica,
Universit\`a degli Studi di Bari}\\
{\small via E. Orabona, 4,  I-70125 Bari, Italy,
 {\tt dambros@dm.uniba.it}} \\ \\
 }
\date{September 19, 2012}

\begin{document}

\maketitle

\begin{abstract} 
  We study a class of systems of 
  quasilinear differential inequalities associated to weakly coercive differential operators and  power reaction terms.
  The main model cases are given by  the $p$-Laplacian operator as well as
  the mean curvature operator in non parametric form.
We prove that if the exponents
  lie under a certain curve, then the system has only the trivial solution.
These results hold without any restriction provided the possible solutions are more regular.
The underlying framework is  the classical Euclidean case as well as
  the Carnot groups setting. 

\medskip

\noindent{\bf Keywords} Quasilinear elliptic systems, a priori estimates, 
 Liouville theorems, Carnot groups. 

\noindent{\bf Mathematics Subject Classification (2010)} 35B40, 35B45, 35J47, 35J92, 35M30, 35R03
\end{abstract}

\section{Introduction}
Liouville theorems for elliptic systems is a classical subject in the theory of partial differential equations.
In recent years, many general results have been obtained by several authors.

In this paper we shall consider a canonical class of elliptic systems of Hamiltonian type of the form,
 \be \label{eq:hamil1}
  \begin{cases} L_1(u)\decl\diver \A_1(x,u,\nabla u)= H_v(x,u,v) \quad on\quad  \RN, \\ \\
 L_2(v)\decl\diver \A_2(x,v,\nabla v)=  H_u(x,u,v) \quad on\quad  \RN,\\ \\
 u\ge 0, \quad v\ge 0\quad on\quad  \RN.
  \end{cases}  \ee
 
Here $L_i , i=1,2$, are quasilinear elliptic operators satisfying a weak coercivity assumption (see below for the
precise assumption) and $H : \mathbb{R}^2\to \mathbb{R}$ is a given function controlled from below by a positive
polynomial in the variables $(u,v).$
Notice that if  in (\ref{eq:hamil1}), at least one of the operators $L_i$ is quasilinear, the problem  is not variational. Plainly the same holds when both operators $L_i$ are linear but $L_1 \neq L_2^*$ (here $L_2^*$ denotes the formal adjoint of $L_2$).

As it is well known one of  the central themes for proving Liouville theorems for problem (\ref{eq:hamil1}) is to find good (possibly sharp) a priori estimates on the weak solutions. These estimates are very difficult to prove  under the weak coercivity assumption of  $L_i , i=1,2$ (see the next section for the precise definition). This is  due mainly to the lack of a weak  Harnack's inequality. 

Results  on a priori estimates of weak solutions  for non coercive quasilinear elliptic systems were first obtained by Mitidieri and Pohozaev in \cite{mit-poh01}, \cite{mit-poh01b} and by Bidaut-V\'eron and Pohozaev in \cite{bid01}
under weak structure assumptions on the differential operators. In those works, particular important examples of systems involving the $p$-Laplacian operator were considered.

Very recently  
these problems have received a renewed interest  for coercive systems.  See for instance Bidaut-V\'eron, Garcia-Huidobro and Yarur  \cite{bid} for a recent interesting contribution in this direction. 

In this paper we shall  focus our attention on weakly coercive elliptic systems of inequalities  of the canonical  form,
\be \label{eq:hamil}
  \begin{cases} \diver \A_1(x,u,\nabla u)\ge v^p \quad on\quad  \RN, \\ \\
\diver \A_2(x,v,\nabla v)\ge u^q \quad on\quad  \RN,\\ \\
 u\ge 0, \quad v\ge 0\quad on\quad  \RN.
  \end{cases}  \ee
Many authors studied system related to (\ref{eq:hamil}). 
We refer the interested reader to Yarur \cite{yar96}, Lair and Wood \cite{lair-wood}, 
Teramoto \cite{ter}, C\^irstea and R\v adulescu \cite{cir-radu}, Yang \cite{yang}
and the references therein. 
The systems studied by these authors 
contain nonautonomous nonlinearities and 
the differential operators are the Laplace operator or the $p$-Laplacian operator. Moreover, many of these papers deal with radial solutions of systems of equations.
Our  results hold even for some classes of nonautonomous systems without assuming the radial symmetry of the possible solutions.
Throughout this paper, in order to avoid cumbersome notations, we  shall focus our attention on (\ref{eq:hamil}). 
However, see Section \ref{sec:na} for some outcomes in the nonautonomous case.

We point out that an  obvious motivation to consider (\ref{eq:hamil}) relies on the fact  that all results on a priori bounds concerning  (\ref{eq:hamil}) apply to the  general systems (\ref{eq:hamil1}) provided,
$$H_v(x,u,v)\ge v^p,$$
$$H_u(x,u,v)\ge u^q.$$
It is  known that if instead of dealing with  system (\ref{eq:hamil}), we study the  scalar model inequality
\be
\diver \A(x,u,\nabla u)\ge u^q \quad on\quad  \RN, \\ \\
u\ge 0, \quad on\quad  \RN, \ee
then, under suitable assumptions on the differential operator  and a natural 
super-homo\-geneity hypothesis on $q$, the problem admits only the trivial 
weak solution.

These kind of Liouville theorems were first proved in  \cite{mit-poh01} under strong assumptions on the differential operators and on the regularity of the possible weak solutions. 

Recently, these results have been extended in different directions  to a wide class
of differential operators and nonlinearities in  Farina and Serrin \cite{FAS}, \cite{FASII} 
 for  $\sigma$-regular distribution solutions which are  {\it locally bounded} and may {\it change sign}.

The {\it locally boundedness} assumption on the solutions has been recently removed in D'Ambrosio and Mitidieri \cite{dam-mit:kato}.

The  non coercive counterpart of (\ref{eq:hamil}), that is
\be \label{eq:hamil-}
  \begin{cases} -L(u) \ge v^p \quad on\quad  \RN, \\ \\
    -L(v)\ge u^q \quad on\quad  \RN,\\ \\
 u\ge 0, \quad v\ge 0\quad on\quad  \RN,
  \end{cases}  \ee
has also been widely studied.  
Mitidieri in \cite{miti2} studies the solutions of (\ref{eq:hamil-})
with $L=\Delta$, proving the existence of the so-called critical Serrin's curve,
\begin{equation*}\max\left\{\frac{p+1}{pq-1}, \frac{q+1}{pq-1} \right\}= \frac{N-2}{2}.\eqno{(S)}\label{hyperbola-}
    \end{equation*}
If the parameters  $(p,q)$ satisfy 
\begin{equation*}\max\left\{\frac{p+1}{pq-1}, \frac{q+1}{pq-1} \right\}\ge \frac{N-2}{2},
    \end{equation*}
 then
the solutions are trivial.

 Dealing with the system of equations (instead of inequalities),
the same kind of result has been proved for radial solutions
in \cite{miti1} by showing the existence of another {\em sharp} critical
Hardy-Littlewood-Sobolev  curve,
\begin{equation*}
\frac 1 {p+1}+\frac 1 {q+1}=\frac  {N-2}N.\eqno{(HLS)}
    \end{equation*}    
    Some contributions on  the so called (HLS) conjecture, see Caristi, D'Ambrosio and Mitidieri \cite{car-dam-mit},  for non radial solutions
 has been obtained by Serrin and Zou \cite{sz}.  
    See also Busca and Manasevich \cite{bm}, 
    Polacik, Quittner and Souplet \cite{PQ}.
    For  ground breaking results on this conjecture see   Souplet \cite{sou}.
    Related results have been recently obtained by
    Phan \cite{phan}.

Later on, problem (\ref{eq:hamil-})  has been studied  in \cite{mit-poh00} 
for  more general operators. Further recent interesting contributions on (\ref{eq:hamil-}) 
and related generalizations can be found in Filippucci \cite{fili}, \cite{fili1},
 and for more general nonlinearities, in D'Ambrosio and Mitidieri \cite{dam-mitST} and  Colasuonno, D'Ambrosio and Mitidieri \cite{col-dam-mit}.
Problems of the type (\ref{eq:hamil-}) has been studied also
in  \cite{car-dam-mit} for linear operators $L$ of higher order as well as 
pseudodifferential operators. 
In \cite{car-dam-mit} the authors use and develop a technique
introduced in \cite{dam-mit-poh} and based on  integral representation 
of the solutions. 
Again, they find critical curves for systems of integral equations and inequalities similar to (S) and (HLS).
These kind of systems are related to the double weighted Hardy-Littlewood-Sobolev inequality studied by Stein and Weiss \cite{stein-weiss}.

In this paper we use a technique  essentially rooted on the classical  test functions method and on some ideas introduced in \cite{dam-mit:kato}.
\medskip

The main goal of this paper is to study weakly coercive systems of the form  (\ref{eq:hamil}) in the framework of Carnot groups for $\sigma$-regular solutions
(see Definition \ref{def:sol} below)
which are not necessarily {\it locally bounded}. 
See the next section for essential preliminaries.
We refer the interested reader to \cite{dam09} and the references therein,
for results related to the scalar case on Carnot groups. 
To the best of our knowledge,  studies concerning coercive systems of type
 (\ref{eq:hamil}) in Carnot groups setting are not present in the literature. 
For the non coercive system see \cite{dam-mitST}.

\medskip
As a sample of our results  in the Euclidean case we have  the following.
\begin{theorem}\label{teo:main0}  Let $\A_i$ be \WPC\ with $p_i>1$ and  $q_i>p_i-1$ ($i=1,2$).
  Let $(u,v)$ be a solution of
\be 
  \begin{cases} \diver \A_1(x,u,\nabla u)\ge v^{q_2} \quad on\quad  \RN, \\ \\
 \diver \A_2(x,v,\nabla v)\ge u^{q_1} \quad on\quad  \RN,\\ \\
 u\ge 0, \quad v\ge 0\quad on\quad  \RN.
  \end{cases}  \label{dis:main}
\ee

 If
\be \max\left\{ q_1\frac{p_2q_2+p_1(p_2-1)}{q_1q_2-(p_1-1)(p_2-1)}+p_1,\ 
  q_2\frac{p_1q_1+p_2(p_1-1)}{q_1q_2-(p_1-1)(p_2-1)}+p_2\right\}\ge N,
  \label{hyperb0}\ee
 then $u\equiv v\equiv 0$\footnote{Throughout this paper we will write
    $u\equiv 0$ for $u=0$ a.e. on $\RN$. }.
\end{theorem} 

As a  special case of  the above theorem we have,
\begin{corollary} Let $\A_i$ be \WPC\ with $p_i>1$ and  $q_i>p_i-1$ ($i=1,2$).
 Let $(u,v)$ be a weak solution of  (\ref{dis:main}). 
  \begin{enumerate}
  \item  If $p_1=p_2=p$, then     $u\equiv v\equiv 0$ provided
$$ \max\left\{q_1\frac{q_2+p-1}{q_1q_2-(p-1)^2}, q_2\frac{q_1+p-1}{q_1q_2-(p-1)^2} \right\}\ge \frac{N-p}{p}. $$
    In particular if $N\le 2p$ then   $u\equiv v\equiv 0$ 
    for any $q_1,q_2>p-1$.
  \item  If $p_1=p_2=p$ and $q_1=q_2=q$, then     $u\equiv v\equiv 0$ provided
    $ q(N-2p)\le (N-p)(p-1)$.
  \item  If $p_1=p_2=2$ and $N\ge 5$, then     $u\equiv v\equiv 0$ provided
    $\max\{q_1, q_2\} \ge q_1q_2\frac{N-4}{2}-\frac{N-2}{2} $.
  \item If $p_1=p_2=2$ and $N\le 4$, then     $u\equiv v\equiv 0$.

    In particular if $L_1=L_2$ is the mean curvature operator and $N\le 4$, 
    then   $u\equiv v\equiv 0$.
  \end{enumerate}
\end{corollary}

By requiring more regularity on one of the operators $L_1$ or $L_2$, we can deduce
Liouville theorems under the only assumption that the exponents satisfy $q_i>p_i-1.$
\begin{theorem} Let $\A_i$ be \WPC\ with $p_i>1$ and  $q_i>p_i-1$ ($i=1,2$).
 Let $(u,v)$ be a weak solution of  (\ref{dis:main}).
 \begin{enumerate}
 \item Assume that $\A_1$ (or $\A_2$) is \SPC.   Then   $u\equiv v\equiv 0$.
   \item Assume that $u$ and   $\A_1$ are radial (or $v$ and $\A_2$).
     Then   $u\equiv v\equiv 0$.
 \end{enumerate}
\end{theorem}

It   is worth  noting  that  in this  paper  we shall perform our  analysis for nonnegative solutions of systems of the form (\ref{dis:main}). This is due to the fact that, if we consider  a slightly more general system  involving nonlinearities that change sign, one cannot expect to reduce the original problem via Kato's inequality (see \cite{dam-mit:kato}) for details) to the study of  solutions with non negative components  as the following simple result shows.
\begin{theorem}\label{teo:diag} Let $(u,v)$ be a weak solution of
   \be 
  \begin{cases} \Delta_p u= \abs v^{q-1}v \quad on\quad  \RN, \\ \\
 \Delta_p v = \abs u^{q-1}u \quad on\quad  \RN,
  \end{cases}  \label{eq:sysanti}\ee
  with $q>1$ and $2\ge p>1$. Then $u=-v$ and $u$ and $v$ are 
  solution of the scalar equation
  $$ -\Delta_p u= \abs u^{q-1}u   \quad on\quad  \RN.  $$

  In particular (\ref{eq:sysanti}) does not admit nonnegative nontrivial  solutions. 
\end{theorem}
From the above result one can infer an analogous result for system
like (\ref{eq:sysanti}) where the operator $\Delta_p$ is replaced by 
$-\Delta_p$.
  See Corollary \ref{cor:sysanti} below.\\

\begin{theorem}\label{teo:diagcurv} Let $(u,v)$ be a weak solution of
   \be 
  \begin{cases} \diver \left( \frac{\nabla u}{\sqrt{1+\abs{\nabla u}^2}}\right)= \abs v^{q-1}v \quad on\quad  \RN, \\ \\
   \diver \left( \frac{\nabla v}{\sqrt{1+\abs{\nabla v}^2}}\right)= \abs u^{q-1}u \quad on\quad  \RN,
  \end{cases}  \label{eq:sysanticurv}\ee
  with $q>1$. Then $u=-v$ and $u$ and $v$ are 
  solution of the scalar equation
  $$ -\diver \left( \frac{\nabla u}{\sqrt{1+\abs{\nabla u}^2}}\right)= \abs u^{q-1}u   \quad on\quad  \RN.  $$

  In particular (\ref{eq:sysanticurv}) does not admit nonnegative nontrivial  solutions. 
\end{theorem}

\bigskip

This paper is organized as follows. In the next section we introduce some notations and preliminaries.
In Section \ref{sec:stime} we prove the main estimates on the solutions of (\ref{dis:main}). The related Liouville
theorems are studied in Section \ref{sec:lio}. Some results for nonautonomous
systems are in Section \ref{sec:na}.
  Section \ref{sec:comp} deals with the study of changing sign solutions of
  some special case of coercive systems.
We end the paper with some open questions.

\section{Notations and  definitions}\label{sec:prehi}

In this paper $\nabla$ and $\abs\cdot$ stand  respectively for the 
usual gradient in $\RN$ and the Euclidean norm.
 $\Omega\subset \RN$ denotes an open set.

In this paper $\nabla$ and $|\cdot|$ stand  respectively for the
usual gradient in $\RN$
and the Euclidean norm.

Let $\mu\in \C(\R^N;\RL)$ be a matrix 
$\mu\decl(\mu_{ij})$, $i=1,\dots,l$, $j=1,\dots,N$.
For $i=1,\dots,l$,  let $X_i$ and its formal adjoint $X_i^*$
be defined as
\begin{equation} X_i\decl\sum_{j=1}^N \mu_{ij}(\xi)\frac{\partial}{\partial \xi_j},
\qquad X_i^*\decl-\sum_{j=1}^N \frac{\partial}{\partial\xi_j}\left(\mu_{ij}(\xi)\cdot\right),
   \label{mu}\end{equation}
and let $\grl$ be the vector field defined by
$\grl\decl (X_1,\dots,X_l)^T=\mu\nabla$
and $\grl^*\decl(X_1^*,\dots,X_l^*)^T$.

For any vector field $h=(h_1,\dots,h_l)^T\in\Cuno(\Omega,\RL)$, we shall use
the following notation $ \diverl(h)\decl\diver(\mu^T h)$,
that is
\begin{equation}  \label{eq:defdiv}
   \diverl(h)=-\sum_{i=1}^l X_i^*h_i=-\grl^*\cdot h.\end{equation}

An assumption that we shall made (which actually is an assumption on the matrix
$ \mu$) is that the operator
$$  \Delta_Gu =\diverl(\grl u)$$
is a canonical sub-Laplacian on a Carnot group 
(see below for a more precise meaning).
The reader, which is not acquainted with these structures, 
can think to the special case 
of $\mu =I$, the identity matrix in $\RN$, that is the usual Laplace operator
in Euclidean setting. Our results are new even in the Euclidean case.

\bigskip

We quote some facts on Carnot groups and refer the interested reader to  
the books of Bonfiglioli, Lanconelli and Uguzzoni \cite{bon-lan-ugu07}, 
Folland and Stein \cite{fol-ste82} and to Folland \cite{fol75}
  for  more detailed information on this subject.

A Carnot group is a connected, simply connected, nilpotent Lie
group $\G$ of dimension $N$ with graded Lie algebra ${\cal
G}=V_1\oplus \dots \oplus V_r$ such that $[V_1,V_i]=V_{i+1}$ for
$i=1\dots r-1$ and $[V_1,V_r]=0$. Such an integer $r$ is called the
\emph{step} of the group.
 We set $l=n_1=\dim V_1$, $n_2=\dim V_2,\dots,n_r=\dim V_r$.
A  Carnot group $\G$ of dimension $N$ can be identified, up to an
isomorphism, with the structure of a \emph{homogeneous Carnot
Group} $(\RN,\circ,\delta_R)$ defined as follows; we
identify $\G$ with $\RN$ endowed with a Lie group law $\circ$. We
consider $\RN$ split in $r$ subspaces
$\RN=\R^{n_1}\times\R^{n_2}\times\cdots\times\R^{n_r}$ with
$n_1+n_2+\cdots+n_r=N$ and $\xi=(\xi^{(1)},\dots,\xi^{(r)})$ with
$\xi^{(i)}\in\R^{n_i}$. 
We shall assume that for any $R>0$ the dilation
$\delta_R(\xi)=(R\xi^{(1)},R^2 \xi^{(2)},\dots,R^r \xi^{(r)})$ 
is a Lie group automorphism.
The Lie algebra of
left-invariant vector fields on $(\RN,\circ)$ is $\cal G$. For
$i=1,\dots,n_1=l $ let $X_i$ be the unique vector field in $\cal
G$ that coincides with $\partial/\partial\xi^{(1)}_i$ at the
origin. We require that the Lie algebra generated by
$X_1,\dots,X_{l}$ is the whole $\cal G$.

We denote with $\grl$ the vector field $\grl\decl(X_1,\dots,X_l)^T$
and we call it \emph{horizontal vector field} and by  $\diverl$ the formal 
adjoint on $\grl$, that is (\ref{eq:defdiv}).
Moreover, the vector
fields $X_1,\dots,X_{l}$ are homogeneous of degree 1 with respect
to $\delta_R$ and in this case
$Q= \sum_{i=1}^r i\,n_i=  \sum_{i=1}^r i\,\mathrm{dim}V_i$ is called the 
\emph{homogeneous dimension} of $\G$.
The \emph{canonical sub-Laplacian} on $\G$ is the
second order differential operator defined by
$$\Delta_G=\sum_{i=1}^{l} X_i^2=\diverl(\grl\cdot )$$ and for $p>1$ the
$p$-sub-Laplacian operator is
$$\Delta_{G,p} u\decl \sum_{i=1}^{l} X_i(\abs{\grl u}^{p-2}X_iu)=
\diverl(\abs{\grl u}^{p-2}\grl u).$$
Since $X_1,\dots,X_{l}$
generate the whole $\cal G$, the sub-Laplacian $\Delta_G$ satisfies the
H\"ormander hypoellipticity  condition.

A nonnegative continuous function $\nu:\RN\to\R_+$ is called a 
\emph{homogeneous norm} on { $\G$}, if 
$\nu(\xi^{-1})=\nu(\xi)$, $\nu(\xi)=0$ if and only if $\xi=0$, and it is
homogeneous of degree 1 with respect to $\delta_R$ (i.e.
$\nu(\delta_R(\xi))=R \nu(\xi)$).
A homogeneous norm $\nu$ defines on $\G$ a \emph{pseudo-distance} defined as
$d(\xi,\eta)\decl \nu(\xi^{-1}\eta)$, which
in general is not a distance.
If $\nu$ and $\tilde \nu$ are two homogeneous norms, then they are equivalent,
that is, there exists a constant
$C>0$ such that $C^{-1}\nu(\xi)\le \tilde \nu(\xi)\le C\nu(\xi)$.
Let $\nu$ be a homogeneous norm, then there exists a constant
$C>0$ such that $C^{-1}\abs\xi\le \nu(\xi)\le C\abs\xi^{1/r}$,
for $\nu(\xi)\le1$.
An example of homogeneous norm is
$  \nu(\xi)\decl\left(\sum_{i=1}^r\abs{\xi_i}^{2r!/i}\right)^{1/2r!}.$

Notice that if $\nu$ is a homogeneous norm differentiable a.e., 
then $\abs{\grl \nu}$ is homogeneous of degree 0 with respect to 
$\delta_R$; hence $\abs{\grl \nu}$ is bounded.

We notice that in a Carnot group, the Haar measure coincides with the Lebesgue measure.

\medskip

Special examples of Carnot groups are the
  Euclidean spaces $\R^Q$.
  Moreover, if $Q\le 3$ then any Carnot group is the ordinary Euclidean
  space $\R^Q$.

 The simplest nontrivial example of a Carnot group
  is the Heisenberg group $\hei^1=\R^3$.
  For an integer $n\ge1$, the Heisenberg group $\hei^n$ is defined as follows:
  let $\xi=(\xi^{(1)},\xi^{(2)})$ with
  $\xi^{(1)}\decl(x_1,\dots,x_n,y_1,\dots,y_n)$ and $\xi^{(2)}\decl t$.
  We endow $\R^{2n+1}$ with  the group law
$\hat\xi\circ\tilde\xi\decl(\hat x+\tilde x,\hat y+\tilde y,\hat t+ \tilde t+2\sum_{i=1}^n(\tilde x_i\hat y_i-\hat x_i \tilde y_i)).$
We consider the vector fields
\[X_i\decl\frac{\partial}{\partial x_i}+2y_i\frac{\partial}{\partial t},\
        Y_i\decl\frac{\partial}{\partial y_i}-2x_i\frac{\partial}{\partial t},
        \qquad\mathrm{for\ } i=1,\dots,n, \]
and the associated  Heisenberg gradient
$ \grh\decl (X_1,\dots,X_n,Y_1,\dots,Y_n)^T$.
The Kohn Laplacian $\lh$ is then the operator defined by
$\lh\decl\sum_{i=1}^nX_i^2+Y_i^2.$
The family of dilations is given by
$\delta_R(\xi)\decl (R x,R y,R^2 t)$ with homogeneous dimension
$Q=2n+2$.
In ${\hei^n}$ a canonical homogeneous norm is defined as
$\abs{\xi}_H\decl \left(\left(\sum_{i=1}^n x_i^2+y_i^2\right)^2+t^2\right)^{1/4}.$

\bigskip

In what follows we shall fix a homogeneous norm $\nu$ and for $R>0$,
 $B_R$ stands for the ball of radius $R>0$ generated by $\nu$, that is 
$B_R\decl \{x : \nu(x)<R\}$ and $A_R$ is
the annulus
$B_{2R}\setminus\overline{B_R}$.
Therefore, by using the dilation $\delta_R$ and the fact that the Jacobian
of $\delta_R$ is $R^Q$, we have
\[ |B_R|=\int_{B_R} dx=R^Q\int_{B_1}dx=w_\nu R^Q\quad
\mathrm{and}\quad|A_R|=w_\nu(2^Q-1)R^Q,\]
where $w_\nu$ is the Lebesgue measure of the unit ball $B_1$ in $\mathbb R^N$. 

\bigskip

Throughout the paper $p\ge1$ and we denote by $W^{1,p}_{L,{loc}}$ the space
$$W^{1,p}_{L,{loc}}(\Omega):=\left\{u\in L^p_{{loc}}(\Omega)\,:\, |\grl u|\in L^p_{{loc}}(\Omega)\right\}.$$

\medskip
Clearly, if $\grl =\nabla$, that is the Carnot group is the trivial
Euclidean space $\RN$, then $ W^{1,p}_{L,{loc}}(\Omega)$
coincides with the usual Sobolev space $W^{1,p}_{{loc}}(\Omega)$.

\bigskip

 In what follows we shall assume that
$\A:\RN\times\R\times\RL\to\RL$ is a Caratheodory
function, that is
for each $t\in\R$ and $\xi\in\RL$ the function $\A(\cdot,t,\xi)$ is measurable; and
for a.e. $x\in\RN$, $\A(x,\cdot,\cdot)$ is continuous.

We consider operators $L$ ``generated'' by $\A$, that is
\[ L(u)(x)=\diverl \left(\A(x,u(x),\grl u(x))\right). \]
Our model cases are the $p$-Laplacian operator, the mean curvature operator and
some related generalizations. See Examples \ref{examp:oper} below.

\begin{definition} \label{ptype} Let $\A:\RN\times\R\times\RL\to\RL$ be a 
  Caratheodory function.
  The function $\A$ is called \emph{weakly elliptic} if it generates a weakly elliptic operator $L$ i.e.
\[
\begin{array}{c}
  \A(x,t,\xi)\cdot \xi\ge 0\ \  \mathit{for\ each\ }x\in\RN,\,t\in\R,\,\xi\in\RL,\\ \\
\A(x,0,\xi)=0\ \ \mathit{or}\ \  \A(x,t,0)=0
\end{array}
\eqno{(WE)}\]

Let $p\ge1$, the function $\A$ is called \WPC\ (weakly-$p$-coercive),
 if $\A$ is (WE) and it generates a weakly-p-coercive operator $L$, i.e.
 if there exists a constant
$k_2>0$ such that
\[ (\A(x,t,\xi)\cdot \xi)^{p-1}\ge k_2^{p-1} \abs{\A(x,t,\xi)}^{p}\ \mathrm{for\ each\ }x\in\RN,\,t\in\R,\,\xi\in\RL.\eqno{(\WPC)}\]

Let $p>1$, the function $\A$ is called {\bf S-$p$-C} (strongly-$p$-coercive),
if there exist $k_1,k_2 >0$ constants such that
\[ (\A(x,t,\xi)\cdot \xi)\ge k_1 \abs \xi^p\ge k_2 \abs{\A(x,t,\xi)}^{p'} \ \mathrm{for\ each\ }x\in\RN,\,t\in\R,\,\xi\in\RL.\eqno{(\SPC)}\]
\end{definition}

For additional information on (WE), \WPC\ and \SPC\ operators see 
\cite{s, bid01, mit-poh01b, dam-mit:kato} and the references
therein.

From now on, if not otherwise specified, $\A$, $\A_1$ and $\A_2$ stands
for \WPC\ function with indices $p$, $p_1$ and $p_2$ respectively.

\begin{definition} \label{def:sol} Let $\Omega\subset\R^N$ be an open set and
 let $f:\Omega\times\R\times\RL \to\R$ be a 
Caratheodory function.  Let $p\ge 1$.  We say that
  $u\in W^{1,p}_{L,loc}(\Omega)$ is a \emph{weak solution}  of
  $$ \diverl \left(\A(x,u,\grl u)\right)\ge f(x,u,\grl u)\qquad on\ \ \Omega,
   $$
 if   $\A(\cdot,u,\grl u)\in L^{p'}_{loc}(\Omega)$,
 $f(\cdot,u,\grl u)\in L^1_{loc}(\Omega)$, 
  and for any nonnegative $\phi\in\C^1_0(\Omega)$ we have
\[ - \int_{\Omega} \A(x,u,\grl u) \cdot\grl\phi \ge\int_{\Omega} f(x,u,\grl u)\phi. \]
\end{definition}
In the literature this kind of solutions are called $\sigma$-regular solutions, see for instance \cite{PS1}.

\begin{definition}
  Let $p_1,p_2\ge 1$ and let $f,g:\Omega\times\R \to\R$ be Caratheodory functions.
  We say that $(u,v)\in W^{1,p_1}_{L,loc}(\Omega)\times W^{1,p_2}_{L,loc}(\Omega)$
  is a \emph{weak solution} of
 \begin{equation}
 \left\{ \begin{array}{ll}
	\diverl(\A_1(x,u,\grl u))\ge { f(x,v)} &
			\qquad\mathrm{on\ }\Omega, \cr \cr
	\diverl(\A_2(x,v,\grl v))	\ge { g(x,u) }&
			\qquad\mathrm{on\ }\Omega.  
			\end{array}\right.  \end{equation}
  if $\A_i(\cdot,u,\grl u)\in L^{p_i'}_{loc}(\Omega)$ ($i=1,2$),
  $f(x,v,\grl v,u,\grl u), g(x,v,\grl v,u,\grl u) \in L^1_{loc}(\Omega)$, 
  and for any nonnegative $\phi\in\C^1_0(\Omega)$ we have
  \begin{equation}\left\{ \begin{array}{l}\displaystyle
    - \int_{\Omega} \A_1(x,u,\grl u) \cdot\grl\phi 
       \ge\int_{\Omega} f(x,v)\phi, \cr\cr\displaystyle
     - \int_{\Omega} \A_2(x,v,\grl v) \cdot\grl\phi 
        \ge\int_{\Omega} g(x,u)\phi.
      \end{array}\right.  \end{equation}  

\end{definition}

\begin{example}\label{examp:oper}
\begin{enumerate}
\item Let $p>1$. The  sub-$p$-Laplacian  operator defined on suitable functions $u$ by,
$$\Delta_{G,p} u=\diverl\left(\abs{\grl u}^{p-2}\grl u\right)$$ 
is an operator generated by
$\A(x,t,\xi):=\abs \xi^{p-2}\xi$ which is  \SPC.

\item  If $\A$ is of mean curvature type, that is
$\A$ can be written as $\A(x,t,\xi):=A(\abs \xi) \xi$ with
  $A:\R\to \R$  a positive bounded continuous function
(see \cite{mit-poh01,bid01}),
then $\A$ is  {\bf W-$2$-C}.

\item The mean curvature operator in non parametric form 
$$Tu: =\diver\left(\frac{\nabla u}{\sqrt{1+\abs{\nabla u}^2}}\right),$$
is generated by $\A(x,t,\xi):=\frac{\xi}{\sqrt{1+\abs{\xi}^2}}$. In this case $\A$ is
  \WPC\ with $1\le p\le2$ and of mean curvature type but it is not  {\bf S-$2$-C}.
\item Let $m>1$. The operator
$$ T_mu\decl\diver\left(\frac{\abs{\gr u}^{m-2}\gr u}{\sqrt{1+\abs{\gr u}^m}}\right)
$$
is \WPC\ for $m\ge p\ge m/2$.

\item Let $p>1$ and define 
$$ L u\decl\sum_{i=1}^N\partial_i\left(\abs{\partial_i u}^{p-2}\partial_iu\right). 
$$
The operator $L$ is \SPC.
\end{enumerate}
\end{example}

Observe that if $\A$ is {\bf S-$p$-C} and $\grl$ is the Euclidean
gradient $\nabla$, then the following important Serrin-Trudinger  inequality holds \cite{s,t,tw}. 
For a proof in the Euclidean setting see, for instance, Maly and Ziemer \cite{MZ},
in the case of Carnot groups see Capogna, Danielli and Garofalo \cite{cdg}.

\begin{lemma}(Weak Harnack Inequality)\label{harnack}
  Let $\A$ be \SPC. 
  If $u\in W^{1,p}_{L,loc}(\RN)$ is a nonnegative weak solution of  
  \begin{equation} {\diverl} (\A(x,u,\grl u))\geq 0, \label{q31}\end{equation} then for any  $\sigma>0$
there exists $c_{H}>0$ independent of $u$ such that for any $R>0$ we have
\begin{equation}
\left (\frac{1}{\abs{B_R}}\int_{B_R}u^\sigma dx\right)^{\frac {1} {\sigma}}\geq\,\, c_{H}\,\,{\rm esssup}_{B_{\frac R 2}} u.\label{q4}
\end{equation}
\end{lemma}
Motivated by the above result, we introduce the following.
\begin{definition}[WH]\label{def:wh}
  We say that the \emph{weak Harnack inequality with exponent $\sigma>0$} 
  holds for $\A$ if for any weak solution $u$ of 
  \be \diverl(\A(x,u,\grl u))\ge 0,\quad u\ge 0,\quad \mathrm{on}\ \Omega\ee
  and any $R>0$ such that $B_{2R}\subset \Omega$ we have
  $$
    \left (\frac{1}{\abs{B_R}}\int_{B_R}u^\sigma dx\right)^{\frac {1} {\sigma}}\geq\,\, c_{H}\,\,{\rm esssup}_{B_{\frac R 2}} u, \eqno{(WH)}
  $$
  with $c_H>0$ independent of $u$ and $R$.
\end{definition}

\section{A priori universal estimates}\label{sec:stime}
In this Section, if not otherwise specified,  $\A, \A_1$ and  $\A_2$ are 
\WPC\ with indices $p>1$, $p_1>1$ and $p_2>1$ respectively.

The following lemma  is a slight variation of a result proved in \cite{dam-mit:kato}. 
For easy reference and  for reader's convenience  we shall include the detailed proof.
\\
\begin{lemma}\label{lem:1sti}
  Let $g\in L^1_{loc}(\Omega)$ be nonnegative  and 
  let $u\in W^{1,p}_{L,loc}(\Omega)$ be a weak solution of
  \be \diverl\left(\A(x,u,\grl u)\right)\ge g,\qquad u\ge0,\quad\ on\ \ \Omega. \label{dis:genloc}\ee
Let $s>0.$ If $u^{s+p-1}\in L^1_{loc}(\Omega)$, then
  \be g u^s,\ \  \A(x,u,\grl u)\cdot \grl u\ u^{s-1}\in L^1_{loc}(\Omega)\label{reg+}\ee
 and for any nonnegative $\phi\in\Cuno_0(\Omega)$ we have,
 \be  \int_\Omega g u^s \phi+c_1 s \int_\Omega\A(x,u,\grl u)\cdot \grl u\ u^{s-1}\phi \le
     c_2 s^{1-p} \int_\Omega u^{s+p-1} \frac{\abs{\grl\phi}^p}{\phi^{p-1}},
      \label{dis:ests}\ee
  where $c_1= 1-\frac{\epsilon^{p'}}{p'k_2}>0$, $c_2=\frac{p^p}{p\epsilon^p}$ and $\epsilon>0 $ is sufficiently small.
\end{lemma}
\begin{remark} i) Notice that from the above result it follows that 
  if  $u\in W^{1,p}_{L,loc}(\Omega)$ is a weak solution of
  (\ref{dis:genloc}), then $g\,u\in L^1_{loc}(\Omega)$.

 ii) The above lemma still holds if we replace the function 
  $g\in  L^1_{loc}(\Omega)$ with a regular Borel measure on $\Omega$.
\end{remark}
\begin{remark}\label{rem:zeros}
  From Lemma \ref{lem:1sti} we deduce that if $p>1$ and 
  $u$ is a weak solution of (\ref{dis:genloc}) and 
  $u$ vanishes on a set $O\subset\Omega$ of positive measure,  then
  $\A(x,u,\grl u)$ must vanishes, modulo sets of measure zero, on $O.$ 
  Indeed, we can choose $0<s<1$ so that
  $u^{s-1}$ is infinity on $O$. Since (\ref{reg+}) holds then necessarily
  $\A(x,u,\grl u)\cdot \grl u=0$ on $O$. Now using the fact that $\A$
  is \WPC\ with $p>1$, we obtain the claim.
\end{remark}

\bp Let $\gamma\in\Cuno(\R)$ be a bounded nonnegative function with
  bounded nonnegative first derivative and let $\phi\in\Cuno_0(\Omega)$
  be a nonnegative test function.

  Applying Lemma 2.2, of \cite{dam-mit:kato} it follows that, 
  \begin{eqnarray*}\lefteqn{
  \int_\Omega g\gamma(u) \phi+ \int_\Omega\A(x,u,\grl u)\cdot\ \grl u\ \gamma'(u)\phi\le - \int_\Omega{\A(x,u,\grl u)}\cdot{\grl\phi}\ \gamma(u)}\\
    &&\le \int_\Omega\abs{\A(x,u,\grl u)}\,\,\abs{\grl\phi} \gamma(u)\\
   &&\le\left(\int_\Omega\abs{\A(x,u,\grl u)}^{p'}\gamma'(u)\phi\right)^{1/p'}
    \left(\int_\Omega\frac{\gamma(u)^p}{\gamma'(u)^{p-1}}  \frac{\abs{\grl\phi}^p}{\phi^{p-1}} \right)^{1/p}\\
 &&\le \frac{\epsilon^{p'}}{p'k_2} \int_\Omega\A(x,u,\grl u)\cdot \grl u\ \gamma'(u)\phi + 
   \frac{1}{p\epsilon^p}\int_\Omega\frac{\gamma(u)^p}{\gamma'(u)^{p-1}} \frac{\abs{\grl\phi}^p}{\phi^{p-1}},
  \end{eqnarray*}
  where $\epsilon>0.$ Notice that all integrals are well defined provided
  $ \frac{\gamma(u)^p}{\gamma'(u)^{p-1}} \in L^1_{loc}(\Omega)$. 
  With a suitable choice of $\epsilon>0,$ for any 
  nonnegative $\phi\in\Cuno_0(\Omega)$  and 
   $\gamma\in\Cuno(\R)$
 such that $\frac{\gamma(u)^p}{\gamma'(u)^{p-1}} \in L^1_{loc}(\Omega),$ 
we have,
  \be  \int_\Omega g\gamma(u) \phi+c_1 \int_\Omega\A\ \grl u \gamma'(u)\phi \le
     \frac{1}{p\epsilon^p} \int_\Omega\frac{\gamma(u)^p}{\gamma'(u)^{p-1}}  \frac{\abs{\grl\phi}^p}{\phi^{p-1}}.
      \label{primastima}\ee

  Now for $s> 0$, $1>\delta>0$ and $n\ge 1$,  define
  \be \gamma_n(t)\decl
  \begin{cases}
    (t+\delta)^s&if\ 0\le t<n-\delta,\\ \\
    \displaystyle cn^s-\frac{s}{\beta-1}{n^{\beta+s-1}}{(t+\delta)^{1-\beta}}&if\ t\ge n-\delta,\\
  \end{cases}  \label{gn}\ee
  where  $c\decl\frac{\beta-1+s}{\beta-1}$ and $\beta>1$ will be chosen later.
  Clearly $\gamma_n\in \C^1$,
  $$\gamma'_n(t)=
  \begin{cases}
    s(t+\delta)^{s-1}&if\ 0\le t<n-\delta,\\ \\
    s{n^{\beta+s-1}}{(t+\delta)^{-\beta}}&if\ t\ge n-\delta,\\
  \end{cases}  $$
  and  $\gamma_n$, $\gamma'_n$ are nonnegative and bounded  with 
  $\norm{\gamma_n}_\infty=cn^s$ and $\norm{\gamma_n'}_\infty=sn^{s-1}$. Moreover 
$$\frac{\gamma_n(t)^p}{\gamma_n'(t)^{p-1}}=
\begin{cases}
  s^{1-p}(t+\delta)^{s+p-1}&\ \  for\ t<n-\delta,\  \\ \\
  \theta(t,n) &\ \ for \ t\ge n-\delta,
\end{cases}
$$
where 
$$\theta(t,n)\decl\frac{(cn^s-\frac{s}{\beta-1}{n^{\beta+s-1}}{(t+\delta)^{1-\beta}})^p}{(s{n^{\beta+s-1}}{(t+\delta)^{-\beta}})^{p-1}}
 \le {(cn^s)^p}{s^{1-p}} n^{-(\beta+s-1)(p-1)}\, (t+\delta)^{\beta(p-1)}.$$
Choosing $\beta\decl \frac{s+p-1}{p-1}$ we have $c=p$, and 
$$\theta(t,n)\le p^p s^{1-p} n^{s p-(\beta+s-1)(p-1)}  (t+\delta)^{s+p-1}= p^p s^{1-p} (t+\delta)^{s+p-1}.$$
Therefore,  for $t\ge 0$ we have,
$$\frac{\gamma_n(t)^p}{\gamma_n'(t)^{p-1}}\le p^p s^{1-p} (t+\delta)^{s+p-1}.$$

Since by assumption  $u^{s+p-1}\in L^1_{loc}(\Omega)$, from (\ref{primastima})
 with $\gamma=\gamma_n$, it follows that
  $$\int_\Omega g\gamma_n(u) \phi+c_1 \int_\Omega\A(x,u,\grl u)\cdot \grl u\ \gamma_n'(u)\phi \le
     \frac{p^p s^{1-p}}{p\epsilon^p}  \int_\Omega (u+\delta)^{s+p-1}\frac{\abs{\grl\phi}^p}{\phi^{p-1}}.$$

Noticing that  $\gamma_n(t)\to (t+\delta)^s$ and $\gamma_n'(t)\to s (t+\delta)^{s-1}$ as 
  $n\to +\infty,$   $g\ge 0$ and $\A\cdot\grl u\ge 0$,  by Beppo Levi theorem 
  we obtain
  $$  \int_\Omega g\, (u+\delta)^s \phi+c_1 s \int_\Omega\A(x,u,\grl u)\cdot \grl u\ (u+\delta)^{s-1}\phi \le
      c_2 s^{1-p}\int_\Omega (u+\delta)^{s+p-1} \frac{\abs{\grl\phi}^p}{\phi^{p-1}},$$
which, by letting $\delta\to 0$,
 completes the proof.
\ep

\begin{remark}\label{rem:smin1}
  The assumption $u^{s+p-1}\in L^1_{loc}(\Omega)$, 
  is not needed for the validity of the statement (\ref{reg+}). 
  Indeed what  really matters  is the assumption  $u^{s+p-1}\in L^1_{loc}(S).$ Here $S$ is the support of $\grl \phi$.
  This observation will be useful when dealing with inequalities on unbounded set.
\end{remark}
\begin{lemma}\label{lem:sti2}
  Let $g\in L^1_{loc}(\Omega)$ be nonnegative  and 
  let $u\in W^{1,p}_{L,loc}(\Omega)$ be a weak solution of (\ref{dis:genloc}).
Let $q>p-1.$ If $u^{q}\in L^1_{loc}(\Omega)$, then
  \be g u^{q-p+1},\ \  \A(x,u,\grl u)\cdot \grl u\ u^{q-p}\in L^1_{loc}(\Omega)\label{reg+2}\ee
 and for any $\varphi\in\Cuno_0(\Omega)$ such that $0\le\varphi\le 1$, we have,
\be\int_\Omega g\varphi^\sigma\le c_3 
  \left(\frac{1}{\abs S}\int_S u^{q}\varphi^\sigma\ \right)^{\frac{p-1}{q}} 
    \left(\frac{1}{\abs S}\int_S \abs{\grl \varphi}^\sigma \right)^{\frac{p}{\sigma}} 
    \abs S,\label{glu}\ee
  where $S$ is the support of $\grl \varphi$,
  $c_3\decl \frac{\sigma^p}{s^{p-1}}\left(\frac{c_2}{c_1k_2}\right)^{1/p'}$
  with
  $\sigma\ge \frac{pq}{q-p+1-s}$, $0<s<\min\{1,q-p+1\}$
  and $c_1, c_2$ as in the above Lemma \ref{lem:1sti}.
\end{lemma}
\bp Let $s>0$ be such that $q\ge s+p-1$. 
  From Lemma \ref{lem:1sti}, for any nonnegative $\phi\in\Cuno_0(\omega)$
  we have  
   \be  \int_\Omega g u^s \phi+c_1 s \int_\Omega\A(x,u,\grl u)\cdot \grl u\ u^{s-1}\phi \le
     c_2 s^{1-p} \int_\Omega u^{s+p-1} \frac{\abs{\grl\phi}^{p}}{\phi^{p-1}}.
     \label{est2}\ee

  Now, let $0<s <\min\{ 1, q-p+1\}$.
  By definition of weak solution and 
   H\"older's inequality with exponent $p'$,
  taking into account that $\A$ is \WPC\ and from  (\ref{est2}) we get,
  \begin{eqnarray} \label{est3}
    \int_\Omega g\phi&\le&\int_S \abs{\A(x,u,\grl u)}\abs{\grl \phi}=
      \int_S\abs{\A}u^{\frac{s-1}{p'}}\phi^{\frac{1}{p'}}\  
      \abs{\grl \phi}u^{\frac{1-s}{p'}}\phi^{-\frac{1}{p'}}
\\
  &\le& \frac{1}{k_2^{1/p'}} \left(\int_S\A(x,u,\grl u)\cdot\grl u\, u^{s-1} \phi\right)^{1/p'}
    \left(\int_S u^{(1-s)(p-1)}\frac{\abs{\grl\phi}^{p}}{\phi^{p-1}} \right)^{1/p}\label{tec13}\\
  &\le&  \frac{1}{k_2^{1/p'}}\left(\frac{c_2}{c_1s^p}\right)^{1/p'}\left(\int_S u^{s+p-1} \frac{\abs{\grl\phi}^{p}}{\phi^{p-1}}\right)^{1/p'} \left(\int_S u^{(1-s)(p-1)}\frac{\abs{\grl\phi}^{p}}{\phi^{p-1}} \right)^{1/p}.\label{est4}
     \end{eqnarray}
  Since  $q>s+p-1$ and $q>p-1$,  applying H\"older
  inequality to (\ref{est4}) with exponents
  $\chi\decl \frac{q}{s+p-1}$ and 
  $y\decl\frac{q}{(1-s)(p-1)}$, we obtain 
  \be
    \int_\Omega g\phi \le c_3' 
\left(\int_S u^{q}\phi \right)^{\delta} 
    \left(\int_S   \frac{\abs{\grl\phi}^{p\chi'}}{\phi^{p\chi'-1}}\right)^{\frac{1}{p'\chi'}} 
\left(\int_S   \frac{\abs{\grl\phi}^{py'}}{\phi^{py'-1}}\right)^{\frac{1}{py'}}, \label{est5pre}\ee 
  where 
  $$\delta\decl \frac{1}{\chi p'}+\frac{1}{yp}=\frac{p-1}{q},\qquad
   c_3'\decl \left(\frac{c_2}{k_2c_1s^p}\right)^{1/p'}.$$
  Next for $\sigma\ge p\chi'$ (notice that, since $ p\chi'> py'$ 
 we have  $\sigma> py'$ )
 we choose $\phi\decl \varphi^\sigma$
  with $\varphi\in\Cuno_0(\Omega)$  such that $0\le\varphi\le 1$.
  Setting $S\decl support(\varphi)$, from (\ref{est5pre}) it follows that
    \be \int_\Omega g \varphi^\sigma \le c_3' \sigma^p 
   \left(\int_S u^{q}\varphi^\sigma \right)^{\delta} 
   \left(\frac{1}{\abs S}\int_S  \abs{\grl\varphi}^{\sigma} 
        \right)^{\frac{p}{\sigma}}  \abs S^{1-\delta},\label{est5}\ee
completing the proof of (\ref{glu}).
\ep

\begin{lemma}\label{lem:apriorigen} Let $q_1>p_1-1$ and $q_2>p_2-1.$
  For any $\sigma>0$ large enough, there exists a constant 
  $c=c(\sigma,q_1,q_2,p_1,p_2,\A_1,\A_2)>0$ such that if $(u,v)$ 
  is weak solution of
  \begin{equation}
    \left\{ \begin{array}{ll}
	\diverl(\A_1(x,u,\grl u))\ge { v^{q_2}} &
			\qquad\mathrm{on\ }\Omega, \cr\\
			
	\diverl(\A_2(x,v,\grl v))	\ge { u^{q_1} }&
			\qquad\mathrm{on\ }\Omega, \\ \\
           v\ge0,\ u\ge 0,
			\end{array}\right.  \label{dis:maingen}  
  \end{equation}
  then for any nonnegative $\varphi\in \Cuno_0(\Omega)$ such that 
  $\norm \varphi_\infty\le 1$, we have
  \begin{eqnarray}
   \int_\Omega u^{q_1}\varphi^\sigma\le c 
   \left(\frac{1}{\abs S}\int_S v^{q_2}\varphi^\sigma\ \right)^{\frac{p_2-1}{q_2}} 
    \left(\frac{1}{\abs S}\int_S \abs{\grl \varphi}^\sigma \right)^{\frac{p_2}{\sigma}} 
    \abs S\label{ulv}\\
 \int_\Omega v^{q_2}\varphi^\sigma\le c 
  \left(\frac{1}{\abs S}\int_S u^{q_1}\varphi^\sigma\ \right)^{\frac{p_1-1}{q_1}} 
    \left(\frac{1}{\abs S}\int_S \abs{\grl \varphi}^\sigma \right)^{\frac{p_1}{\sigma}} 
    \abs S\label{vlu}\\
 \int_\Omega u^{q_1}\varphi^\sigma\le c\ \abs S\ 
   \left(\frac{1}{\abs S}\int_S \abs{\grl \varphi}^\sigma \right)^{\frac{1}{\sigma}\,\frac{p_1\delta_2+p_2}{1-\delta_1\delta_2}} \label{estu}\\
\int_\Omega v^{q_2}\varphi^\sigma\le c\ \abs S\ 
   \left(\frac{1}{\abs S}\int_S \abs{\grl \varphi}^\sigma \right)^{\frac{1}{\sigma}\,\frac{p_2\delta_1+p_1}{1-\delta_1\delta_2}},  \label{estv}
  \end{eqnarray}
  where $S\decl support( \varphi)$,
  $\delta_1\decl \frac{p_1-1}{q_1}$ and $\delta_2\decl \frac{p_2-1}{q_2}$.
\end{lemma}

\bp From Lemma \ref{lem:sti2} we immediately obtain (\ref{ulv}) and (\ref{vlu}).

  Since $$\int_S v^{q_2}\varphi^\sigma\le \int_\Omega v^{q_2}\varphi^\sigma,$$ by inserting (\ref{vlu}) in
  (\ref{ulv}) we get (\ref{estu}). In a similar way we can prove (\ref{estv}).
\ep

\begin{lemma}\label{lem:estR} Let $q_1>p_1-1$ and $q_2>p_2-1.$ 
  There exists $c>0$ such that if  $(u,v)$ is weak solution of 
  (\ref{dis:maingen}), $x\in\Omega$ and $B_{2R}(x)\subset\subset\Omega$, then
  \be \left( \mint_{B_R(x)} v^{q_2}\right)^{1/q_2}\le c 
    R^{-\frac{p_1q_1+p_2(p_1-1)}{q_1q_2-(p_1-1)(p_2-1)}},\label{estvR}\ee
  \be \left( \mint_{B_R(x)} u^{q_1}\right)^{1/q_1}\le c 
    R^{-\frac{p_2q_2+p_1(p_2-1)}{q_1q_2-(p_1-1)(p_2-1)}},\label{estuR}\ee
\begin{equation}
  \mint_{B_R(x)} v^{q_2} \le c \left(\mint_{B_{2R}(x)} u^{q_1} 
    \right)^{\frac{p_1-1}{q_1}}  R^{-p_1},\qquad  \mint_{B_R(x)} u^{q_1} \le c 
    \left(\mint_{B_{2R}(x)} v^{q_2} \right)^{\frac{p_2-1}{q_2}}  R^{-p_2}.
\end{equation}
\end{lemma}
\bp Let $\phi_0\in\Cuno_0(\R)$ be such that
  $$\phi_0(t)=0\ for\ \abs t\ge 2,\   \phi_0(t)=1\ for \ \abs t\le 1,\ 
  and \ 0\le\phi_0\le 1.$$
  For $R>0$, we define $\phi_R(y)\decl\phi_0(\frac{\nu(x^{-1}y)}{R})$.

  By choosing $\varphi=\phi_R$ in Lemma \ref{lem:apriorigen},
  and observing that
  $$\left(\frac{1}{\abs S}\int_S \abs{\grl \varphi}^\sigma \right)^{\frac{1}{\sigma}}= 
  \left(\frac{1}{\abs{B_{2R}(x)\setminus B_R(x)}}\int_{B_{2R}(x)\setminus B_R(x)} \abs{\grl \phi_R(y)dy}^\sigma \right)^{\frac{1}{\sigma}}
  = R^{-1}c(\phi_0,Q,\nu,\sigma),$$
the claim follows from (\ref{ulv}), (\ref{vlu}), (\ref{estu}) and  (\ref{estv}).
\ep

An immediate consequence of the above results are the  following universal estimates
on the solutions of (\ref{dis:maingen}).
\begin{lemma}\label{lem:univdist}
  Let $q_1>p_1-1$ and $q_2>p_2-1.$ 
  There exists $c>0$ such that for any $(u,v)$ weak solution of 
  (\ref{dis:maingen}), $x\in\Omega$ and $R=dist(x,\partial\Omega)/2$,
  we have
  \be \left( \mint_{B_R(x)} v^{q_2}\right)^{1/q_2}\le c 
    \ dist(x,\partial \Omega)^{-\frac{p_1q_1+p_2(p_1-1)}{q_1q_2-(p_1-1)(p_2-1)}},\label{est:v}\ee
  $$ \left( \mint_{B_R(x)} u^{q_1}\right)^{1/q_1}\le c \ dist(x,\partial \Omega)^{-\frac{p_2q_2+p_1(p_2-1)}{q_1q_2-(p_1-1)(p_2-1)}},$$
  \begin{equation}\label{dis:acc}
  \mint_{B_R(x)} v^{q_2} \le c \left(\mint_{B_{2R}(x)\setminus B_R(x)} u^{q_1} \right)^{\frac{p_1-1}{q_1}}  dist(x,\partial \Omega)^{-p_1},\ee
\be   \mint_{B_R(x)} u^{q_1} \le c \left(\mint_{B_{2R}(x)\setminus B_R(x)} v^{q_2} \right)^{\frac{p_2-1}{q_2}}  dist(x,\partial \Omega)^{-p_2}.
\end{equation}

  Moreover if (WH) holds for  $\A_1$ (resp. $\A_2$),
  then for a.e. $x\in\Omega$ we have
  \begin{eqnarray}
      u(x) \le c \ dist(x,\partial \Omega)^{-\frac{p_2q_2+p_1(p_2-1)}{q_1q_2-(p_1-1)(p_2-1)}}\label{estux}\\
  \left(resp.\ \  v(x) \le c \ dist(x,\partial \Omega)^{-\frac{p_1q_1+p_2(p_1-1)}{q_1q_2-(p_1-1)(p_2-1)}}\right).\label{estvx}
  \end{eqnarray}
\end{lemma}
\begin{remark}
  The estimates proved in the above lemma above are sharp. 
  That is,  estimates (\ref{estux}) and (\ref{estvx}) cannot be improved.
To see this consider  $\Omega\decl]0,+\infty[\times \R^{N-1}$.
  In this case  $dist(x,\partial\Omega)=x_1$.
  Let $q_1>p_1-1>0$ $q_2>p_2-1>0$.  It is easy to see that the functions $u, v$ defined by
  $$u(x)\decl x_1^{-\frac{p_2q_2+p_1(p_2-1)}{q_1q_2-(p_1-1)(p_2-1)}}\quad and \quad
v(x)\decl x_1^{-\frac{p_1q_1+p_2(p_1-1)}{q_1q_2-(p_1-1)(p_2-1)}}$$
  solve  the system, 
$$
 \left\{ \begin{array}{ll}
       \Delta_{p_1}u= \lambda { v^{q_2}} &
			\qquad\mathrm{on\ }\Omega, \cr
       \Delta_{p_2}v= \mu { u^{q_1} }&
			\qquad\mathrm{on\ }\Omega.  
			\end{array}\right.
$$
where $\lambda,\mu>0$ are suitable positive constants.
\end{remark}

\begin{remark}The above results, 
  Lemma \ref{lem:1sti}, \ref{lem:sti2}, \ref{lem:apriorigen}
  \ref{lem:estR}, \ref{lem:univdist}, continue to hold (with the same proof as above) even in the case 
  $p_1=1$ or $p_2=1.$   
\end{remark}
\section{Liouville Theorems}\label{sec:lio}
In this Section, if not otherwise specified,  $\A, \A_1$ and  $\A_2$ are 
\WPC\ with indices $p>1$, $p_1>1$ and $p_2>1$ respectively.\\

Our main result in this section is the following.
\begin{theorem}\label{teo:main} 
  Let $(u,v)$ be a weak solution of 
\be 
  \begin{cases} \diverl \A_1(x,u,\grl u)\ge v^{q_2} \quad on\quad  \RN, \\ \\
 \diverl \A_2(x,v,\grl v)\ge u^{q_1} \quad on\quad  \RN,\\ \\
 u\ge 0, \quad v\ge 0\quad on\quad  \RN.
  \end{cases}  \label{dis:mainG}
\ee
 If
\be \max\left\{ q_1\frac{p_2q_2+p_1(p_2-1)}{q_1q_2-(p_1-1)(p_2-1)}+p_1,\ 
  q_2\frac{p_1q_1+p_2(p_1-1)}{q_1q_2-(p_1-1)(p_2-1)}+p_2\right\}\ge Q,
  \label{hyperb}\ee
 then $u\equiv v\equiv 0$.

If $p_1=p_2=p$, then  (\ref{hyperb}) becomes
$$ \max\left\{q_1\frac{q_2+p-1}{q_1q_2-(p-1)^2}, q_2\frac{q_1+p-1}{q_1q_2-(p-1)^2} \right\}\ge \frac{Q-p}{p},
$$
which in turn, when $q_1=q_2=q$ becomes
$$ q(Q-2p)\le (Q-p)(p-1).$$
\end{theorem} 
\begin{remark}
  \begin{enumerate}
  \item  Notice that (\ref{hyperb})  can be rewritten as
\be \max\left\{ p_1(\frac1\delta_1-1),p_2(\frac1\delta_2-1)\right\}\ge
   \frac{Q-p_1-p_2}{\delta_1\delta_2}-Q, \label{newhyperb}\ee
where $\delta_i\decl \frac{p_i-1}{q_i}$, $i=1,2$.

\item  From (\ref{newhyperb}) we easily see that if $Q\le p_1+p_2$ then there are
  no nontrivial solution of (\ref{dis:mainG}) for any $q_1$ and $q_2$ such that
  $q_i>p_i-1$, $i=1,2$.
  \end{enumerate}
\end{remark}

\noindent{\bf Proof of Theorem \ref{teo:main}.} 
   Define
$$f_0\decl Q-q_1\frac{p_2q_2+p_1(p_2-1)}{q_1q_2-(p_1-1)(p_2-1)} \ \ \ 
and
 \ \ \ g_0\decl 
Q-q_2\frac{p_1q_1+p_2(p_1-1)}{q_1q_2-(p_1-1)(p_2-1)}.$$

Applying Lemma \ref{lem:1sti} to the first inequality of (\ref{dis:mainG})
with $s\decl q_1-p_1+1$ we have
 \be\label{tec43}  \int_{\rn} v^{q_2} u^{q_1-p_1+1} \phi+c_1 s \int_{\rn}\A_1(x,u,\grl u)\cdot \grl u\ u^{q_1-p_1}\phi \le
     c_2 s^{1-p_1} \int_{\rn} u^{q_1} \frac{\abs{\grl\phi}^{p_1}}{\phi^{p_1-1}}.\ee
Taking $\phi=\phi_R$ a standard cut off function as in the proof of 
Lemma \ref{lem:estR} and using  the above estimates (\ref{estuR}) we have
\begin{equation}  \label{estuq}
  \int_{\rn} u^{q_1} \frac{\abs{\grl\phi}^{p_1}}{\phi^{p_1-1}}\le c  R^{f_0-p_1},
\end{equation}
which, together with (\ref{tec43}) implies
\be  \int_{B_R} v^{q_2} u^{q_1-p_1+1} +c_1 s \int_{B_R}\A_1(x,u,\grl u)\cdot \grl u\ u^{q_1-p_1} 
     \le c R^{f_0-p_1}.
 \label{est12}\ee
 Similarly we obtain that
 \be  \int_{B_R} u^{q_1} v^{q_2-p_2+1} +c\int_{B_R}\A_2(x,v,\grl v)\cdot \grl v\ v^{q_2-p_2}
     \le c R^{g_0-p_2}.
 \label{est13}\ee
  Since the hypothesis (\ref{hyperb}) can be written as
  $\min\{f_0-p_1,g_0-p_2\}\le 0$, we deduce that at least 
 one of the exponent in  (\ref{est12}) and (\ref{est13}) is nonpositive.
 Without loss of generality we assume that 
  $f_0\le p_1$.
  Then, from (\ref{est12}) we obtain that 
  \begin{equation}
    \label{AL1}  \A_1(x,u,\grl u)\cdot \grl u\ u^{q_1-p_1}\in L^1(\RN).    
  \end{equation}

  Now arguing as in Lemma \ref{lem:1sti}, using $u^{q_1-p_1+1} \phi$ as test function, we have
  \begin{eqnarray}
     \lefteqn{\int_{\rn} v^{q_2} u^{q_1-p_1+1} \phi+(q_1-p_1+1) \int_{\rn}\A_1(x,u,\grl u)\cdot \grl u\ u^{q_1-p_1}\phi \le} \\
   &&  \le \int_S u^{q_1-p_1+1}\abs{\A_1}\abs{\grl \phi}\\ 
   && \le
     \left(\int_S\A_1(x,u,\grl u)\cdot \grl u\ u^{q_1-p_1}\phi \right)^{1/p_1'}
     \left(\int_S u^{q_1} \frac{\abs{\grl\phi}^p_1}{\phi^{p_1-1}}\right)^{1/p_1}\\ 
   && \le   \left(\int_S\A_1(x,u,\grl u)\cdot \grl u\ u^{q_1-p_1}\phi \right)^{1/p_1'} M, \label{tec33}
  \end{eqnarray}
  where $S$ is the support of $\grl \phi$ and in the last inequality 
  we have used (\ref{estuq}) and the fact that  $f_0-p_1\le 0$.
  Choosing again $\phi=\phi_R$ as above, from (\ref{tec33}) we obtain
  $$\int_{B_R}\A_1(x,u,\grl u)\cdot \grl u\ u^{q_1-p_1} \le
  c \left(\int_{B_{2R}\setminus B_R}\A_1(x,u,\grl u)\cdot \grl u\ u^{q_1-p_1} \right)^{1/p_1'},$$
  which, with the information (\ref{AL1}), by letting $R\to +\infty$, implies that
  $$ \A_1(x,u,\grl u)\cdot \grl u\ u^{q_1-p_1} \equiv 0
   \ \ \ on\ \RN.$$

  From Remark \ref{rem:zeros} we have that $\A_1(x,u,\grl u)\equiv 0$, 
  and consequently, by the first inequality of (\ref{dis:mainG}), 
  $v\equiv 0$. The same conclusion holds for $u$. 
\ep

\begin{theorem}\label{teo:wh} Assume that either $\A_1$ or $\A_2$ satisfy the
  weak Harnack inequality  (WH).
  If $q_i>p_i-1$, $i=1,2$, and $(u,v)$ is a weak solution of (\ref{dis:mainG}),
  then $u\equiv v\equiv0$.
\end{theorem}
\noindent{\bf Proof of Theorem \ref{teo:wh}.} Assume that for $\A_1$ the weak Harnack inequality holds. From estimate (\ref{estuR}) we have
$$ \sup_{B_{R/2}} u\le c R^{-\frac{p_2q_2+p_1(p_2-1)}{q_1q_2-(p_1-1)(p_2-1)}}.$$
By letting $R\to+\infty$ in the above inequality, the claim will follows.
\ep

\begin{corollary}  Assume that either $\A_1$ is  {\bf S}-$p_1$-{\bf C} or
  $\A_2$ is {\bf S}-$p_2$-{\bf C}.
  If $q_i>p_i-1$, $i=1,2$, and $(u,v)$ is a weak solution of (\ref{dis:mainG}),
  then $u\equiv v\equiv0$.
\end{corollary}
\bp Assume that $\A_1$ is {\bf S}-$p_1$-{\bf C}.
  From Lemma \ref{harnack} we have that for 
  the nonnegative solutions of $\diverl(\A(x.u.\grl u))\ge 0$ 
  the weak Harnack inequality holds. Thus Theorem \ref{teo:wh} applies.
\ep

In what follows we need of the following.
\begin{lemma}\label{lem:inf} Let $w\in L^1_{loc}(\RN)$ be nonnegative.
  For $r>0$ define\footnote{We recall that $\nu(x)$ denotes  the norm of $x$.}
  $$ m_w(r)\decl {\rm ess}\inf_{\nu (x)>r} w(x).$$
  Let $f:]0,+\infty[\to ]0,+\infty[$ be a positive continuous
   nondecreasing function such that 
  $f(w)\in L^1_{loc}(\RN)$. 
  If
  $$ \liminf_R \mint_{B_R} f(w(x) dx = 0, $$
  then $\lim_{t\to 0} f(t) = 0$ and  for any $r>0$,
  $ m_w(r)=0.$
\end{lemma}
\noindent{\bf Proof.} 
  Let $r_0>0$ and $R>r_0$. We have
  \begin{eqnarray*}
    \mint_{B_R} f(w)  \ge 
    \frac{1}{\omega_\nu R^Q} \int_{B_R\setminus B_{r_0}} f(w) \ge \frac{1}{ R^Q} (R^Q-r_0^Q) f\left(m_u(r_0)\right).
  \end{eqnarray*}
  Letting $R\to+\infty$ in the last chain of inequalities, we obtain 
  $0\ge f(m_u(r_0))$.
\ep
\begin{corollary}\label{cor:inf} Let $q_i>p_i-1$, $i=1,2$, and 
  let $(u,v)$ be a weak solution of (\ref{dis:mainG}).
  Let $r>0$ and set
  $$ m_u(r)\decl {\rm ess}\inf_{\nu (x)>r} u(x), \qquad m_v(r)\decl{\rm ess}\inf_{\nu (x)>r} v(x).  $$
  Then  for any $r>0$ we have $m_u(r)=m_v(r)=0$.
\end{corollary}
\noindent{\bf Proof.} 
  From estimates  (\ref{estvR}) and (\ref{estuR}) it follows  that
  $$\left( \mint_{B_R} v^{q_2}\right)^{1/q_2}\ \  and \ \ 
  \left( \mint_{B_R} u^{q_1}\right)^{1/q_1}$$ vanish for $R\to +\infty$.
  By choosing $f(t)=t^{q_2}$ and $f(t)=t^{q_1}$ in Lemma \ref{lem:inf}, the claim follows. 
\ep

The next result deals with radial solution in the Euclidean framework.
 We need  the following.
\begin{definition} We say that $\A$ is radial if there exists
  a  Caratheodory function
  $A:\R\times\R\times\R\to\R$ such that 
  for any $u\in \Cuno(\RN)$
  radial, that is $u(x)=u(\abs x)$, we have that
 $$\A(x,u,\gr u) = A(\abs x, u (\abs x), u'(\abs x))\frac{x}{\abs x}.$$
\end{definition}  
\begin{remark}
  If $\A$ is radial and $u=u(\abs x)$ is a radial  solution of
  $$ \diver(\A(x,u,\gr u)\ge f(\abs x) \quad on\ \RN,$$
  then $u$ solves
  $$ (r^{N-1} A(r,u(r),u'(r)))'\ge r^{N-1} f(r), \quad for\ \  r>0.$$
\end{remark}

\begin{theorem}\label{teo:rad} Let $\G$ be the Euclidean group $\RN$.
  Let $q_i>p_i-1$, $i=1,2$, and 
  let $(u,v)$ be a weak solution of (\ref{dis:mainG}).
  Assume that either  
  $\A_1$ is radial and $u$ is a $\Cuno$ radial function
  or  $\A_2$ is radial and $v$ is a $\Cuno$ radial function.
  Then $u\equiv v\equiv0$.
\end{theorem}
\bp Assume that $\A_1$ is radial and $u$ is a $\Cuno$ radial function.
  Therefore $u$ solves
  $$ (r^{N-1} A(r,u(r),u'(r)))'\ge 0 \quad for \  \ r>0.$$
    Integrating between $0$, and $r$ we have
    $r^{N-1} A(r,u(r),u'(r))\ge 0$
  which together with the weakly ellipticity (WE) of $\A$ yields
 $u'(r) \ge 0$ for $r>0$.

  We proceed by contradiction assuming that $u\not \equiv 0$. Hence
  there exists $r_0>0$  such that $u(r_0)>0$. Since $u$ is nondecreasing
  we have that $u(r)\ge u(r_0)>0$ for any $r>r_0$.
  This  contradicts    Corollary \ref{cor:inf}.
\ep

\subsection{Nonautonomous systems}\label{sec:na}
In this Section we briefly show how the ideas developed in the preceding
 Sections can be  
 employed to study  a class of nonautonomous systems.

We consider
\be 
  \begin{cases} \diverl \A_1(x,u,\grl u)\ge H(x)\, v^{q_2} \quad on\quad  \RN, \\ \\
 \diverl \A_2(x,v,\grl v)\ge K(x)\, u^{q_1} \quad on\quad  \RN,\\ \\
 u\ge 0, \quad v\ge 0\quad on\quad  \RN.
  \end{cases}  \label{dis:mainGna}
\ee

We  assume that $H$ and $K$ are measurable functions
such that for a.e. $x\in \RN$
\begin{equation}
  \label{eq:hypHK} 
  H(x)\ge \frac{c_H}{(1+\nu(x))^{\alpha_1}},\qquad
  K(x)\ge \frac{c_K}{(1+\nu(x))^{\alpha_2}},
\end{equation}
where $\alpha_1,\alpha_2\ge 0$ and $c_H,c_K>0$.

We have the following.
\begin{theorem}\label{teo:mainGna} Let $q_i>p_i-1$, $i=1,2$ and let $(u,v)$ be a weak solution of (\ref{dis:mainGna}).
  \begin{enumerate}
  \item 
   If
\begin{eqnarray*} \max\left\{ 
q_1\frac{(p_2-\alpha_2)q_2+(p_1-\alpha_1)(p_2-1)}{q_1q_2-(p_1-1)(p_2-1)}+p_1,\ 
q_1\frac{(p_2-\alpha_2)q_2+(p_1-\alpha_1)(p_2-1)}{q_1q_2-(p_1-1)(p_2-1)}+\alpha_2,
  \right.\nonumber\\ 
\left. q_2\frac{(p_1-\alpha_1)q_1+(p_2-\alpha_2)(p_1-1)}{q_1q_2-(p_1-1)(p_2-1)}+p_2, \   
q_2\frac{(p_1-\alpha_1)q_1+(p_2-\alpha_2)(p_1-1)}{q_1q_2-(p_1-1)(p_2-1)}+\alpha_2\right\}> Q,
  \label{hyperbna}\end{eqnarray*}
 then $u\equiv v\equiv 0$.
\item
  Assume that either,
   $$\A_1 \ is\  {\bf S}-p_1-{\bf C}\quad and \quad 
   (p_2-\alpha_2)q_2+(p_1-\alpha_1)(p_2-1)>0,$$
   or
  $$    \A_2\  is\  {\bf S}-p_2-{\bf C}
   \quad and \quad (p_1-\alpha_1)q_1+(p_2-\alpha_2)(p_1-1)>0 .$$
  Then $u\equiv v\equiv0$.
  \end{enumerate}
\end{theorem}
\begin{remark}
  The above result complements and refines those obtained in \cite{ter}. 
  Indeed, in \cite{ter}
  the author proved existence and nonexistence of radial solutions for
  a special systems of equations involving  $p$ and $q$-Laplacian operators.
\end{remark}

Since the proof of Theorem \ref{teo:mainGna} follows the same ideas as in Theorem \ref{teo:main} and Theorem \ref{teo:wh}, we shall be brief.
We first state some estimates which have an interest in themselves.
These  are the analogues of 
those obtained in Lemma \ref{lem:estR}, which in this case read as
  \be \left( \mint_{B_R(x)} H v^{q_2}\right)^{1/q_2}\le c 
    R^{-\frac{(p_1-\alpha_1)q_1+(p_2-\alpha_2)(p_1-1)}{q_1q_2-(p_1-1)(p_2-1)}-\frac{\alpha_1}{q_2}},\label{estvRna}\ee
  \be \left( \mint_{B_R(x)} Ku^{q_1}\right)^{1/q_1}\le c 
    R^{-\frac{(p_2-\alpha_2)q_2+(p_1-\alpha_)(p_2-1)}{q_1q_2-(p_1-1)(p_2-1)}-\frac{\alpha_2}{q_1}
},\label{estuRna}\ee
for any $x\in\RN$ and $R>0$.
By using the hypotheses on $H$ and $K$, it foolows that
 \be \left( \mint_{B_R(x)} v^{q_2}\right)^{1/q_2}\le c 
    R^{-\frac{(p_1-\alpha_1)q_1+(p_2-\alpha_2)(p_1-1)}{q_1q_2-(p_1-1)(p_2-1)}},\label{estvRna2}\ee
  \be \left( \mint_{B_R(x)} u^{q_1}\right)^{1/q_1}\le c 
    R^{-\frac{(p_2-\alpha_2)q_2+(p_1-\alpha_)(p_2-1)}{q_1q_2-(p_1-1)(p_2-1)}
}.\label{estuRna2}\ee
  Inequalities (\ref{estvRna}), (\ref{estuRna}), (\ref{estvRna2}), and 
  (\ref{estuRna2}),    allow us to argue as in the proofs of
  Theorem \ref{teo:main} and Theorem \ref{teo:wh} obtaining the claim.

  \begin{remark}
Notice that, the same kind of results hold if we replace the pointwise conditions
(\ref{eq:hypHK}) with the weaker ones
\begin{equation*}
  \label{eq:hypHKweak} 
  \left(\mint_{R<\nu(x)<2R}H^{-\sigma}\right)^{1/\sigma} \le C_H R^{\alpha_1},\qquad
 \left(\mint_{R<\nu(x)<2R}K^{-\sigma}\right)^{1/\sigma} \le C_K R^{\alpha_2},
\end{equation*}
for $R$ large and  suitable   $\sigma>0$. 
These kind of integral hypotheses on the coefficients have already been used
in \cite{apriori} and \cite{dam-mit:kato}. We refer the  
interested reader to those works for further details.
  \end{remark}

\section{Systems with changing sign nonlinearities}\label{sec:comp}

Let $L w=\diverl\A(x,w,\grl w)$ be as above. Let $(u,v)$ be a solution of
\begin{equation}
 \left\{ \begin{array}{ll}
	L(u)\ge f(v) &
			\qquad\mathrm{on\ }\Omega, \cr
			\\
	L(v)	\ge g(u)&
			\qquad\mathrm{on\ }\Omega. 
			\end{array}\right.  \label{dis:sysgen}  
\end{equation}
Clearly the pair $(u,v)$ solves the inequality
\be L(u)+L(v)\ge g(u)+f(v) \qquad\mathrm{on\ }\Omega.\label{dis:comp1}\ee
In particular if $L$ is odd (that is $L(-u)=-L(u)$ ), then setting $w\decl -u$
and $\tilde g(t)\decl - g(-t)$,  the couple $(v,w)=(v,-u)$
solves the inequality
\be L(v) - f(v) \ge L(w)-  \tilde g(w) \qquad\mathrm{on\ }\Omega.
\label{dis:comp}\ee 

In order to study (\ref{dis:comp}) we need of the following.
\begin{definition}
  We say that $\A:\RN\times \RL\to \RL$ is monotone if 
  \be  (\A(x,\xi)-\A(x,\eta)) \cdot (\xi-\eta) \ge 0
   \qquad for\ \ x\in\RN,\ \xi,\eta\in\RL. \label{hypAmon}\ee
  We say that $\A$ is \MPC\ (monotone $p$-coercive) if
  $\A$ is monotone and there exists $k>0$ such that
  \begin{equation}  \label{CWPC}
    \left( \A(x,\xi)-\A(x,\eta)\cdot (\xi-\eta) \right)^{p-1}\ge k^{p-1}
    \abs{\A(x,\xi)-\A(x,\eta) }^{p} \qquad for\ \ x\in\RN,\ \xi,\eta\in\RL.
  \end{equation}
\end{definition}

Examples of \MPC\ operators can be found for instance in
 D'Ambrosio, Farina and Mitidieri \cite{dam-far-mit}.
  Among others we  have,
\begin{example}
  \begin{enumerate}
  \item  Let $1<p\le 2,$ then  the function 
  $\A(\xi)\decl \abs {\xi}^{p-2}\xi$ generates an  \MPC\ operator.

  \item The mean curvature operator is \MPC\ for $1<p\le 2$.
  \end{enumerate}
\end{example}

In \cite{dam-far-mit} the authors prove the following comparison result.
\begin{theorem}\label{teo:comp}  Let $\A$ be \MPC\ and $q>\min\{1,p-1\}$. Let (u,v) be a weak solution of
\be \diverl\left(\A(x,v,\grl v)\right)-\abs {v}^{q-1}v
\ge \diverl\left(\A(x,u,\grl u)\right)-\abs{u}^{q-1}u\qquad on\ \RN.
\label{dis:cmp}\ee
Then $v\le u$ a.e. on $\RN$.
\end{theorem}

\begin{remark} Notice that the above claim $v\le u$ a.e. on $\RN,$ is fairly easy to prove under locally bounded assumption of the weak solutions.
For further details on this point see  \cite{dam-far-mit}.
  \end{remark}

As consequence of Theorem \ref{teo:comp},  we prove the following.
\begin{theorem}\label{th:q=p} Let $\A$ be \MPC, $q>\min\{1,p-1\}$ and let $(u,v)$ be
  a weak solution of 
  \begin{equation}
    \left\{ \begin{array}{ll}
	  \diverl\left(\A(x,\grl u)\right)\ge \abs v^{q-1}v &
			\qquad\mathrm{on\ }\RN, \cr
			\\
	 \diverl\left(\A(x,\grl v) \right)\ge\abs u^{q-1}u &
			\qquad\mathrm{on\ }\RN.  
			\end{array}\right.  \label{dis:sysq}  
  \end{equation}
  Then $u+v\le 0$ a.e. on $\RN$.
 
  Moreover, if $\A$ is odd, that is $\A(x,-\xi)=-\A(x,\xi)$, and  $(u,v)$
  solves also the equation in (\ref{dis:sysq}), then $u=-v$ a.e. on $\RN$ and $u$ solves
 $$- \diverl\left(\A(x,\grl u)\right)= \abs u^{q-1}u\qquad \mathrm{on}\ \ \RN.$$
\end{theorem}
\bp Arguing as at the beginning of this section,
  $(-u,v)$ is a solution of (\ref{dis:cmp}). 
  Hence by Theorem \ref{teo:comp} it follows that $v\le -u.$
  This completes the first part of the proof.

Now, if  $(u,v)$ is a solution of (\ref{dis:sysq}) with equality sign, then 
  $(-u,-v)$ solves the same equations. By the first part of this claim  we deduce that $-u-v\le 0$, thereby concluding the proof.
\ep

\begin{corollary}\label{cor:sysanti}
 Let $\A$ be \MPC\ and odd. Let $q>\min\{1,p-1\}$ and let $(u,v)$ be
  a weak solution of 
  \begin{equation}
  \left\{ \begin{array}{ll}
	 - \diverl\left(\A(x,\grl u)\right)= \abs v^{q-1}v &
			\qquad\mathrm{on\ }\RN, \cr
			\\
	- \diverl\left(\A(x,\grl v) \right) = \abs u^{q-1}u &
			\qquad\mathrm{on\ }\RN.  
			\end{array}\right.  \label{dis:sysqanti}  
  \end{equation}
  Then $u\equiv v$ a.e. on $\RN$.
\end{corollary}
\bp The claim follows by observing that $(-u,v)$ solves the system
  (\ref{dis:sysq}) with equality signs. Hence the claim follows from 
  Theorem \ref {th:q=p}. \ep

\begin{corollary}
     Let $\A$ be \MPC, $q>\min\{1,p-1\}$ and let $(u,v)$ be nonnegative
  weak solution of (\ref{dis:sysq}), then $u\equiv v\equiv 0$.
\end{corollary}

\subsection*{Some open questions}
Obviously, there are many interesting questions one may ask on the problems
considered  in this paper.
We find the following particularly intriguing.
\begin{enumerate}
  \item We  proved our main result under the restrictive hypothesis
     $$q_1>p_1-1,$$ and $$q_2>p_2-1.$$ The main question is that  if the same result as in Theorem \ref{teo:main}
holds  under the weaker  natural assumption (see \cite{bid})     $$q_1q_2>(p_1-1)(p_2-1).$$
    Clearly, this problem is meaningful if one knows that Harnack's inequality
    does not hold.
  \item Is the curve in Theorem \ref{teo:main} sharp? 
 We believe that the answer is no.
  \item If the nonlinearities are not power functions, 
    is there any condition {\em a la}
    Keller-Osserman implying that the only solution of the system (\ref{dis:main}) is the trivial one?
  \item Does the Liouville's theorem hold if in system  (\ref{eq:sysanticurv})
    we replace the power functions with increasing nonlinearities?
  \item What about Theorem \ref{teo:diag} when the power of the nonlinearity is not the same?
\end{enumerate}

\subsubsection*{Acknowledgment}
We acknowledge the support of
 MIUR National Research Project: {\it Quasilinear Elliptic Problems and Related Questions}.

\bibliographystyle{amsplain}

\end{document}